\def\eps{\varepsilon}

\def\s{\sigma}

\def \dd#1{{\bf#1}}

\def\cl#1{{\cal#1}}



\def\ouv#1{\smash{\mathop{#1}\limits^{\lower 1pt\hbox
{$\scriptscriptstyle\circ$}}}}

\def\hfl#1#2{\smash{\mathop{\hbox to 12mm{\rightarrowfill}}
\limits^{\scriptstyle#1}_{\scriptstyle#2}}}


\long\def\eno#1#2{\par\smallskip{\bf{#1}}{\it\ {#2}}\par\medskip}

\def\tit#1{\vskip 5mm plus 1mm minus 2mm {\tir #1}
		\vskip 3mm plus 1mm minus 2mm}

\def\stit#1{\vskip 3mm plus 1mm minus 2mm {\bf{#1}}
		\smallskip}

\long\def\dem#1#2{{\bf {#1}}{\ {#2}$\diamondsuit$}\medskip}

\font\tir=cmbx10 at 12pt

\def\ref#1#2#3#4{{\bf #1}{\ #2}{\it ,\ #3}{,\ #4}\medskip}


\def \picture #1 by #2 (#3){\midinsert \centerline 
{\vbox to #2{\hrule width #1 heigth 0pt 
depth 0pt \null \vfill \special {picture #3}}}\endinsert }

\def\scaledpicture #1 by #2 (#3 scaled #4) {{
\dimen0 =#1 \dimen1 =$2
\divide \dimen0 by 1000 \multiply \dimen0 by #4
\divide \dimen1 by 1000 \multiply \dimen1 by #4
\picture \dimen0 by \dimen1 (#3 scaled $4)}}

\def\figure #1 #2 #3 {\midinsert \vglue 3mm 
{\vbox to #3 {\hrule width 6cm height 0cm depth 0cm \vfill
{\special {picture #1 scaled #2}}}}\vglue 2mm \endinsert}

\magnification=1200

\overfullrule=0pt

{\centerline {\tir {The standard deviation effect}}}

{\centerline {\tir {(or why
one should sit  first base playing blackjack)}}}

\bigskip
\bigskip

{\centerline {{\bf E. Mu\~noz-Garc\'\i a $^*$,  R. P\'erez-Marco}\footnote {*} 
{  UCLA, 
Dept. of Mathematics, 
405, Hilgard Ave., Los Angeles, CA-90095-1555, 
USA, e-mail: munoz@math.ucla.edu, ricardo@math.ucla.edu.}}}

\bigskip

{\bf Abstract.} {\it For a balanced cardcounting system we study
the random variable of the true count after a number of cards 
are removed from the remaining deck and 
we prove a close formula for its standard deviation. 
As expected, the formula shows that the standard deviation 
increases with the number of cards removed.
This creates a "standard deviation effect" with a two
fold consequence: 
longer long run and presumably larger fluctuations of the 
bankroll, but 
a small gain in playing accuracy 
for the player sitting third base.
The opposite happens for the player sitting first base. 
Thus  the optimal position
in casino blackjack in terms of shorter long run is first base.
}

\bigskip


\medskip


\bigskip
\bigskip

{\centerline {\bf Contents}}

\bigskip
	{\bf 1) Introduction.}\par
	\medskip

	{\bf 2) Standard deviation of the true count.}\par
	\medskip

		{\parindent=2cm

		\item{}a) Proof of the true count theorem.\par
		\smallskip

		\item{}b) The true count algebra.\par
		\medskip

		\item{}c) The true count theorem without induction.\par
		\medskip
		
		\item{}d) The true count standard deviation formula.\par
		\medskip

	}

	{\bf 3) Long run.}\par
	\medskip

		{\parindent=2cm

		\item{}a) Kelly for binomial games.
		\smallskip
		\item{}b) Kelly for fuzzy advantage.\par
		\medskip
		\item{}c) Long run.\par
		\medskip
		}

	{\bf 4) Practical gambling.}\par
	\medskip

		{\parindent=2cm

		\item{}a) Comparison of first and third base and 
head-on play.\par
		\smallskip
		\item{}a1) Full table.
		\smallskip
		\item{}In company of ploppys
		\smallskip
		\item{}a2) Head-on play.
		\smallskip
		\item{}b) Absolute magnitude of the standard deviation.\par
		\medskip

		}

	{\bf Bibliography.}\par
	\medskip

\null 
\vfill
\eject

{\bf Disclaimer.}
\medskip
Cardcounters are asked to forgive us the simplistic 
presentation of their activity aimed to (non-cardcounter) mathematicians, and 
mathematicians those elementary computations in the last section 
aimed to (non-mathematician) cardcounters.

Not being gamblers or probability specialists this article
is written for fun and with no pretention.

\bigskip
\bigskip

\tit {1) Introduction.}

This article is about casino blackjack and not about   
tournament blackjack.
Cardcounters in the game of blackjack use 
different count systems to keep 
track of the ratio of favorable and unfavorable cards that 
remain on the deck. It is well known now and documented ([Th]) 
that one can get an edge (of the order of 2\% ) 
over the house by increasing the 
bets in favorable counts and keep them to a minimum in 
unfavorable ones. The first person who published and 
analyzed this finding, mathematically and with computer simulations, 
was E. Thorp (who at the time was a graduate
student in the mathematics department in UCLA).

\bigskip

One of the main questions that a professional cardcounter
player faces is the money management. He plays with 
a given bankroll. He has to decide the amount to bet at each 
moment taking into account the limit of the table and the 
upper bound fixed by his own bankroll. Obviously the maximum
bet planned in his strategy should be inferior to the maximum
of the table. But how much to bet depending on how favorable is
the count ?  The goal is to follow a betting pattern that 
minimizes the risk (that is loosing the whole bankroll) but 
maximizes the growth rate.
The well known sharp strategy is Kelly's criterion : The bet
is a fraction of the total bankroll equal
to the advantage you have. This maximizes 
the expected exponential rate of growth of the bankroll.
Kelly's criterion can be proven to be optimal in a very strong 
sense: Any other strategy will have a longer expected time 
to achieve a given amount of the bankroll 
(L. Breiman [Br] \footnote {*}{At the time also
a UCLA professor}). For stablishing this type of results 
one makes the assumption that there is no minimal unit
bet. In section 2 we give a simple derivation of Kelly criterion. 
For a basic bibliography we refer to [RT].
When the advantage is not know precisely but it is 
a random variable, Kelly criterion also applies (see section 3.c): 
One should bet
according to the expected value of the advantage. The 
standard deviation of the exponential rate of growth of the bankroll 
turns out to increase in a significant way with 
the standard deviation of the advantage. This makes 
longer the "long run" (see section 3) and presumably
induces larger fluctuations. Thus one should prefer playing conditions that 
minimize the standard deviation of the advantage.

\bigskip

In practice, count systems associate a value to each card.
We restrict the discussion here to balanced counts, that is, those
for which the total sum of values of a cards in a complete 
deck is zero. In blackjack only the numerical value  
of the card matters. A very simple 
balanced count system called Hi-Lo gives the weight $-1$ to 
high cards (10,J,Q,K,A), the weight $+1$ to low 
cards (2,3,4,5,6) and the weight $0$ to medium cards (7,8,9).
The cardcounter player adds up the weights of the cards
as they are revealed. The numerical value obtained is 
the {\it running count} (RC). It tells the player how favorable is
the remaining deck (from the Hi-Lo values you can tell that 
high cards are in favor of the player, low cards are in 
favor of the dealer). Of course how favorable is the deck 
depends not directly on the RC but on the {\it true
count} (TC) that is the ratio of RC by the number of cards
remaining (one can also divide by the number of decks remaining,
this is just a change of units, and it is what is done in 
practice). The quantity $52.TC$ (i.e. the running count in 
deck units, $52$ is the total number of cards in a deck) 
gives approximatively the percentage in edge gained by 
the player due to the unbalanced composition of the deck.

\bigskip

The advantage of the house with a game with favorable rules
facing a perfect play with no counting  from 
the player (this is called {\it basic strategy}) is of the 
order of $0.5 \%$ (see for example, [Gr2] p.139). 
Thus it is recommended to beat the minimum
when $52.TC < 2$ and to bet $52.TC$ units
when $ 52.TC > 1$.

\bigskip

The blackjack tables are semi-circular. The dealer stands on the
flat side and at most seven players (in some tables five players)
sit around the circular border. 
Cards are dealt from the left side to the right side of 
the dealer. The first position at the left
of the dealer is the {\it first base} and the last position at the 
right of the dealer is the {\it third base}. The dealer deals the 
cards from left to right.

It seems to be a common belief in the blackjack literature that
there is no preferred position for the card counter and 
that the company of other players at the table has no influence 
(see [HC] p.68, 
[Wo] p. 220, 
for example). In this article we prove 
that the positions in the 
blackjack table are not all equivalent. 
There is an advantage 
in terms of shorter long run (and presumably of smaller 
fluctuations of the bankroll) for 
the player to be as close as possible to first base. More precisely,
a player betting according to Kelly's criterium will have a larger 
fluctuations of the exponential rate of growth of his bankroll 
if he sits on first base than if he
sits in third base. The main noticeable effect of this is 
in the "long run". A player sitting third base will need to 
play about $2\% $ more favorable hands than a player sitting 
first base in order to achieve the same standard deviation 
for the exponential rate of growth of his bankroll. 
On the other hand, from the edge point of view, a
player sitting closer to third base has a slight 
playing advantage over the players at his right.

The main goal of this article is to explain theoretically  
these differences in the position of the players.
They both have the same origin: 
the {\it standard deviation effect}.  In a few words and 
in a simplistic form, the 
standard deviation effect appears when one has to take an 
irreversible decision with respect to a future situation. 
More the lag is important, more the decision has chances 
to be incorrect.

\bigskip

The game of blackjack is played as follows. First the players decide
their bets. A card counter will do this according to the value of the true 
count, and he will bet following Kelly's criterium. 
Then the dealer will deal
two cards to each player and two for him. The cards are dealt face up
or down, this is not important for our purposes (in all cases there is
always a hidden card, the so-called hole card of the dealer). Then
the players will make their play decisions (to split, double, hit or stand) 
in turn from first to third base, and play the hand. 
When the moment of the play decision comes
the true count of the player has changed because he has 
seen some of the new cards (at least the two he has received).
Thus he may not be betting in a favorable situation.
The so-called "True count theorem" justifies the action. The expected 
value of the true count after some cards have been removed from 
the remaining deck is the same as the true count before removing 
the cards. As far as we know, the 
proof of this theorem in the blackjack context 
was publicished for the first time by 
Abdul Jalib M'Hall ([JM]) 
in a message posted to the rec.gambling.blackjack
newsgroup.

\eno {Theorem (True count theorem, A. J. M'Hall).}{For any 
balanced count, the expected value of the true count after 
several cards (but not all) are removed  from the remaining 
deck is the value of the true count before removing the 
cards.}

This type of result goes back to the 
origins of probability theory.
There is a well known problem:
Each person in a room is asked to take a ticket from a box. One indicates
the winner. Should you try to take your ticket as soon as possible or 
as late as possible? We assume that everyone waits that all tickets 
have been chosen to look at their ticket (otherwise if they look and 
reveal the result the problem is not the same).
The answer is that it doesn't matter of course. The expected value 
of the probability of being a winner (running count) stays the same 
at all moments.

\medskip

Coming back to the card counting problem, even if the expected 
value of the true count is independent of the number of cards 
removed, the standard deviation of the true count will increase
as we remove cards from the deck. This results in particular 
from the exact formula
proved in the theorem below. Thus the precise knowledge
of the true count "dilutes" as more cards are played.
This makes that often 
the advantage the player has 
when his turn of play comes (and also 
when the dealer's turn comes) differs from the expected value
at the moment of betting.
This induces a longer long run.
Since the standard deviation increases with the number of cards 
played, the player sitting in third base will be systematically
hurt by a larger standard deviation, thus will experience a 
longer long run.
This is one of the consequences of the {\it standard deviation effect}.
On the other hand, when the player in third base has to play, he 
is closer to the dealer's play. Thus, compared to the player in first
base, he knows more accurately the true count at the moment the 
dealer will play. Thus he can adjust 
more effectively his play. 
He has a supplementary bet advantage 
than the 
player in first base because he can deviate in a more efficient
way from basic strategy. This is the second consequence of the 
standard deviation effect. The gain from deviating from basic 
strategy is small.

\null
\vfill
\eject

Now we present the "True count standard deviation formula" which
proves the increasing of the standard deviation with the number of 
cards removed from the remaining deck.
Let $\s_n$ be the standard
deviation of the true count after having removed $n\geq 1$ cards.

\eno {Theorem (True count standard deviation formula).}
{Let $N$ be the number of cards remaining in the deck. After 
removing $1\leq n <N$ cards, the standard deviation of the true count is 
$$
\s_n =\sqrt {N-1 \over N-n} \  \sqrt n  \ \s_1.
$$
} 

The combinatorics in the proof of the true count standard 
deviation formula reveal a beautiful set of identities
that we call the {\it true count algebra}. Each one of 
these identities has a probabilistic interpretation. 
We only develop in this article
the minimum set of identities necessary in the proof of the main 
theorem. The distribution of the true count does depend on the 
composition of the deck and the weights of the counting system.
This family of distributions contains the hypergeometric distribution.
We don't know a reference for the true count distribution in 
the classical literature. The authors will be grateful if any 
reader can provide one.

We have the following corollaries (using also the results in 
section 2).
 
\eno {Corollary 1.}{The standard deviation $\s_n$ is strictly 
increasing with $n$, even faster when the deck contains fewer 
cards.}

\eno {Corollary 2.}{ For a given player, at the moment 
of taking his betting decision,
the standard deviation of the true count at the moment of his play
decision will be larger if he sits further away from the first base.}

\eno {Corollary 3.}{The theoretical long run of players
sitting further away from the first base are larger.}

\eno {Corollary 4.}{In terms of shorter theoretical long run, the 
optimal seat is first base.}

\eno {Corollary 5.}{For a given player, at the moment of taking
his playing decision, the standard deviation of the true count at the
moment of the play of the dealer will be smaller if he sits 
further away from first base.}

\eno {Corollary 6.}{The playing advantage of players 
sitting further away from the first base is increased.}

\eno {Corollary 7.}{In terms of playing advantage the 
optimal seat is third base.}

We mainly discuss the effects 
on long run 
that is the most relevant one. There is also 
an effect on larger fluctuations of the bankroll
and risk of ruin that is more difficult to analyze 
with precision.
Classically 
the trade-off between betting efficiency and playing efficiency
is something the card counter has to consider carefully when 
choosing his system (these quantities differ for different 
systems). Empirical "proper balances" have been proposed
(see [Gr1] p.40-49), but of course the proper balance between 
this two different quantities is up to each player and his goal.

\bigskip

\stit {Practical gambling.}

In the last section  we present a closer study of how the 
playing conditions do influence the standard deviation of the true 
count, thus the long run.  
The true count standard deviation formula only gives the relation
of $\s_n$ with $\s_1$. The actual value of $\s_1$ does depend on 
the count system and on the remaining distribution of weights 
in the remaining cards. One has a quite accurate approximate 
formula only depending on the count system and the number 
of cards on the deck
$$
\s_1 \approx {\Sigma_0 \over N}
$$
where $\Sigma_0$ is the standard deviation of weights used in the 
count system.
We discuss some relevant
consequences for practical play. Some amusing consequence
of the formulas is that for a continuous model of the deck 
the standard deviation goes to $+\infty$ when $N\to 0$.
Fortunately (!) the casinos do not use to practice a $100 \% $
penetration on the decks. Of course, in normal 
conditions an important penetration
will increase in a substantial amount the playing advantage
that will become the predominant effect. Nevertheless
it is not excluded that 
casinos could device a set of rules that look advantageous to 
the players (according to the classical literature) 
but induce large fluctuations that will 
wipe out the bankrolls of card counters in the long run.

\bigskip

The methods exposed in this article provide also information 
about higher moments of the true count. The combinatorics one 
faces is more involved. We plan to study this in the future.

\bigskip

{\bf Acknowledgements.}

The authors are grateful to T. Ferguson, T. Liggett, 
R. Schonmann and B. Rothschild for instructive discussions in 
game theory, probability and combinatorics (domains in 
which the authors are amateurs). The authors thank 
also to F. Cheah and Sonia for their numerical experiments and
gambling experience respectively, and also for motivating
the present article.

\null
\vfill
\eject

\tit {2) Standard deviation of the true count.}

\stit {a) Proof of the true count theorem.}

To set up the notations we recall first Abdul Jalib M'Hall's 
proof [JM]. It is enough to prove the result when we remove one card.

We denote by $N$ the number of cards remaining in the deck. 
We denote by $w$ the different possible weights of the cards 
given by the count system. Let $s_w$ be the total number of cards 
in the full deck with weight $w$.
Since the system count is balanced, we have 
$$
\sum_w w s_w =0.
$$
Let $l_w$ be the total number
of cards in the remaining deck with weight $w$ (so $1=\sum_w l_w /N$). 
Then if we denote 
by $k_w = s_w -l_w$ we have that 
$$
R=\sum_w w k_w =-\sum_w w l_w
$$
is the running count. The true count is $T=R/N$. 

We compute the expected value $T_1$ of the true count after 
removing $1$ card. It is given by
$$
\eqalign {
T_1&=\sum_w {R+w \over N-1} {l_w \over N} \cr
&=\sum_w {R\over N-1} {l_w \over N} +\sum_w {w\over N-1} {l_w \over N} \cr
&={R\over N-1} -{R\over (N-1) N}\cr
&={R\over N}\cr
&=T\cr}
$$
and the result follows.

\stit {b) The true count algebra.}

For $1\leq p \leq N$ we define 
$$
l_w^{w_1 \ldots  w_p} =l_w -|\{ 1\leq i \leq p ; w_i =w \} | ,
$$
that is the number of cards of weight $w$ remaining once
$p$ cards of weight $w_1, \ldots , w_p$ have been removed.
Observe that we have 
$$
\sum_w {l_w^{w_1 \ldots w_p} \over N-p} =1.
$$
These coefficients have a  beautiful and rich combinatorics.
We call this system of identities the true count algebra. 
The name is chosen because most relevant formulas have a 
"true count" probabilistic interpretation.
We concentrate here into the relevant identities in order to 
compute the standard deviation. The full algebra will be studied 
in a forthcoming article.

\eno {Lemma 1.}{For $0 \leq p \leq N-2$ we have
$$
S_0=\sum_w {l_w^{w_1 \ldots w_p} \over N-p} \ 
{l_{v_0}^{w_1 \ldots w_p w} \over N-p-1} ={l_{v_0}^{w_1 \ldots w_p} 
\over N-p}.
$$
}
One can translate this relation 
in a probabilistic statement. 
Picking at random 
one card on the remaining deck, the probability that this card 
has a weight $v_0$ is the same than the probability of this 
same occurrence using the same deck depleted from one card at 
random. Or, in card counters terms, the running count for 
weight $v_0$ does not change when we remove one card at random.

\dem {Proof.}{We split the sum into a sum over all $w\not= v_0$ 
and the term for $w=v_0$. Observe that for $w\not= v_0$ we have
$$
l_{v_0}^{w_1 \ldots w_p w}=l_{v_0}^{w_1 \ldots w_p }.
$$
Also for $w=v_0$ we have 
$$
l_{v_0}^{w_1 \ldots w_p v_0}=l_{v_0}^{w_1 \ldots w_p }-1.
$$
Thus we have
$$
\eqalign {
S_0 &= \left ( \sum_{w\not= v_0} {l_w^{w_1 \ldots w_p} \over N-p} 
\ {l_{v_0}^{w_1 \ldots w_p } \over N-p-1} \right )
+{l_{v_0}^{w_1 \ldots w_p} \over N-p} 
\ {l_{v_0}^{w_1 \ldots w_p v_0} \over N-p-1} \cr
&={l_{v_0}^{w_1 \ldots w_p } \over N-p-1}
\left (\sum_{w\not= v_0} {l_w^{w_1 \ldots w_p} \over N-p}
\right ) +{l_{v_0}^{w_1 \ldots w_p} \over N-p} 
\ {l_{v_0}^{w_1 \ldots w_p }-1 \over N-p-1} \cr
&={l_{v_0}^{w_1 \ldots w_p } \over N-p-1}
\left ( 1- {l_{v_0}^{w_1 \ldots w_p} \over N-p} \right )
+ {l_{v_0}^{w_1 \ldots w_p} \over N-p} 
\ {l_{v_0}^{w_1 \ldots w_p }-1 \over N-p-1} \cr
&=l_{v_0}^{w_1 \ldots w_p } \left ( {1\over N-p-1} -{1\over 
(N-p)(N-p-1)}\right )\cr
&={l_{v_0}^{w_1 \ldots w_p} 
\over N-p}\cr}
$$
q.e.d.}

\eno {Lemma 2.}{Let $p\geq 0$, $q\geq 0$ such that $p+q \leq N-2$.
We have the formula
$$
\eqalign {
S_q &=\sum_w {l_w^{w_1 \ldots w_p} \over N-p} \ 
{l_{v_0}^{w_1 \ldots w_p w} \over N-p-1} \ 
{l_{v_1}^{w_1 \ldots w_p v_0 w} \over N-p-2} \ 
\ldots 
{l_{v_q}^{w_1 \ldots w_p v_0 \ldots v_{q-1} w} \over N-p-q-1} \cr
&={l_{v_0}^{w_1 \ldots w_p} 
\over N-p}\ {l_{v_1}^{w_1 \ldots w_p v_0} 
\over N-p-1}\ \ldots 
{l_{v_q}^{w_1 \ldots w_p v_0 \ldots v_{q-1} } 
\over N-p-q}\cr
}
$$
Clearing denominators one can also write the less cumbersome 
formula
$$
\eqalign {
& \sum_w l_w^{w_1 \ldots w_p} \ 
l_{v_0}^{w_1 \ldots w_p w}  \ 
l_{v_1}^{w_1 \ldots w_p v_0 w} \ 
\ldots 
l_{v_q}^{w_1 \ldots w_p v_0 \ldots v_{q-1} w}  \cr
&=
(N-p-q-1)
\ l_{v_0}^{w_1 \ldots w_p} \ 
l_{v_1}^{w_1 \ldots w_p v_0} \ \ldots 
l_{v_q}^{w_1 \ldots w_p v_0 \ldots v_{q-1} }\cr
}
$$
}

One can translate the identity of the lemma  
in a probabilistic statement. 
Picking at random 
$q$ cards on the remaining deck one after the other, 
the probability that the cards 
have  weights  $v_0$,  $v_1$ ...and $v_q$ in this 
order 
is the same than the probability of this 
same occurrence using the same deck depleted from one card at 
random. Or, in card counters terms, the running count for
an ordered clump of $q+1$ cards to   
have respective weights $v_0$,  $v_1$, ... and $v_q$ 
does not change when we remove one card at random.

\dem {Proof.}{The proof is done by induction on $q\geq 0$.
The result holds for $q=0$ because of lemma 1.
We assume that the result holds for $q\geq 0$ and we prove
it for $q+1$.
We split the sum 
$$
\sum_w l_w^{w_1 \ldots w_p} \ 
l_{v_0}^{w_1 \ldots w_p w}  \ 
l_{v_1}^{w_1 \ldots w_p v_0 w} \ 
\ldots 
l_{v_{q+1}}^{w_1 \ldots w_p v_0 \ldots v_{q} w}  
$$
into the sum for $w\not= v_0$,
$$
A=\sum_{w\not= v_0} l_w^{w_1 \ldots w_p} \ 
l_{v_0}^{w_1 \ldots w_p w}  \ 
l_{v_1}^{w_1 \ldots w_p v_0 w} \ 
\ldots 
l_{v_{q+1}}^{w_1 \ldots w_p v_0 \ldots v_{q} w} 
$$
and the term of the sum for $w=v_0$
$$
B=l^{w_1 \ldots w_p}_{v_0}\  l^{w_1 \ldots w_p v_0}_{v_0} \ 
l^{w_1 \ldots w_p v_0 v_0}_{v_1} \ l^{w_1 \ldots w_p v_0 v_1 v_0}_{v_2}
\ldots l^{w_1 \ldots w_p v_0 \ldots v_q v_0}_{v_{q+1}}
$$
For $w\not=v_0$ we have $l_{v_0}^{w_1 \ldots w_p w}=l_{v_0}^{w_1 \ldots w_p}$,
so
$$
A=l_{v_0}^{w_1 \ldots w_p} 
\sum_{w\not= v_0} l_w^{w_1 \ldots w_p} \   \ 
l_{v_1}^{w_1 \ldots w_p v_0 w} \ 
\ldots 
l_{v_{q+1}}^{w_1 \ldots w_p v_0 \ldots v_{q} w} 
$$
Since for $w\not=v_0$ we have $l_w^{w_1 \ldots w_p}=l_w^{w_1 \ldots w_p v_0}$,
we have 
$$
A=l_{v_0}^{w_1 \ldots w_p} 
\sum_{w\not= v_0} l_w^{w_1 \ldots w_p v_0} \   \ 
l_{v_1}^{w_1 \ldots w_p v_0 w} \ 
\ldots 
l_{v_{q+1}}^{w_1 \ldots w_p v_0 \ldots v_{q} w} 
$$
Using the induction hypothesis for $q$ we conclude that 
$$
\eqalign {
A=&l_{v_0}^{w_1 \ldots w_p}   [ (N-(p+1)-q)\  l_{v_1}^{w_1 \ldots w_p v_0}
\ldots l_{v_{q+1}}^{w_1 \ldots w_p v_0\ldots v_q} -\cr
&-l_{v_0}^{w_1 \ldots w_p v_0} \ 
l^{w_1 \ldots w_p v_0 v_0}_{v_1} \ l^{w_1 \ldots w_p v_0 v_1 v_0}_{v_2}
\ldots l^{w_1 \ldots w_p v_0 \ldots v_q v_0}_{v_{q+1}} ]  \cr
&=(N-p-(q+1)) \ l_{v_0}^{w_1 \ldots w_p} \ 
l_{v_1}^{w_1 \ldots w_p v_0}
\ldots l_{v_{q+1}}^{w_1 \ldots w_p v_0\ldots v_q}
-B \cr}
$$
And the result follows.}

From the probabilistic interpretation we obtain that there 
are corresponding formulas for the removal of $k$ cards 
instead of only one. We prove these formulas  
by induction on $k$ using lemma 1 and lemma 2.

\eno {Lemma 3.}{Let $k\geq 1$ and $0\leq p\leq N-k-1$. We have
$$
\sum_{i_1, \ldots , i_k} {l_{i_1}^{w_1 \dots w_p} \over N-p} 
\ {l_{i_2}^{w_1 \dots w_p i_1} \over N-p-1} \ 
{l_{i_3}^{w_1 \dots w_pi_1 i_2 } \over N-p-2} \ldots 
{l_{i_k}^{w_1 \dots w_p i_1 \ldots i_{k-1}} \over N-p-k+1}\ 
{l_{v_0}^{w_1 \dots w_p i_1 \ldots i_{k}} \over N-p-k} 
={l_{v_0}^{w_1 \dots w_p} \over N-p}
$$
}

\eno {Lemma 4.}{Let $p\geq 0$, $q\geq 0$ and $k\geq 1$
such that $p+k+q \leq N-1$. We have 
$$
\eqalign{
&\sum_{i_1 , \ldots , i_k} {l_{i_1}^{w_1 \ldots w_p} \over N-p}
\ {l_{i_2}^{w_1 \ldots w_p i_1} \over N-p-1} \ldots 
{l_{i_k}^{w_1 \ldots w_p i_1\ldots i_{k-1}} \over N-p-k-1}\ 
{l_{v_0}^{w_1 \ldots w_p i_1\ldots i_{k-1} i_k} \over N-p-k} \cr
& \ \ \ \ \ \ \ \ \ \ \ \ \ \ \ \ \ \ \ \ \ \ \ \ \ 
{l_{v_1}^{w_1 \ldots w_p v_0 i_1\ldots i_{k}} \over N-p-k-1}\ 
\ldots 
{l_{v_q}^{w_1 \ldots w_p v_0 \ldots v_{q-1} i_1\ldots i_{k}} \over N-p-k-q} \cr
&={l_{v_0}^{w_1 \ldots w_p} 
\over N-p}\ {l_{v_1}^{w_1 \ldots w_p v_0} 
\over N-p-1}\ \ldots 
{l_{v_q}^{w_1 \ldots w_p v_0 \ldots v_{q-1} } 
\over N-p-q}\cr
}
$$
}

There are more general formulas, but we don't need them for 
the purpose of this article. Next lemma is an application of 
the formulas.

\eno {Lemma 5.}{Let $W=\{ w \}$ be the set of weights 
and $f: W \to \dd R$ a real valued function. The expected 
value of $f(w)$, for any $1\leq n\leq N$ and 
$1\leq j\leq n$, is given by 
$$
\bar f =\sum_w f(w) \ {l_w \over N} =\sum_{w_1 \ldots w_n}
f(w_j) \ {l_{w_1} \over N} \ {l_{w_2}^{w_1} \over N-1}\ldots
{l_{w_n}^{w_1 \ldots w_{n-1}} \over N-n+1} \ .
$$
}

\dem {Proof.}{We have 
$$
\eqalign {
&\sum_{w_1 \ldots w_n}
f(w_j) {l_{w_1} \over N} \ {l_{w_2}^{w_1} \over N-1}\ldots
{l_{w_n}^{w_1 \ldots w_{n-1}} \over N-n+1} \cr
&=\sum_{w_1 \ldots w_j} f(w_j) {l_{w_1} \over N} 
\ {l_{w_2}^{w_1} \over N-1}\ldots 
{l_{w_j}^{w_1\ldots w_{j-1}} \over N-j+1}
\ \sum_{w_{j+1} , \ldots w_n} {l_{w_{j+1}}^{w_1\ldots w_{j}} \over N-j}
\ldots {l_{w_n}^{w_1\ldots w_{n-1}} \over N-n+1} \cr
&=\sum_{w_1 \ldots w_j} f(w_j) {l_{w_1} \over N} 
\ {l_{w_2}^{w_1} \over N-1}\ldots 
{l_{w_j}^{w_1\ldots w_{j-1}} \over N-j+1} \cr
&=\sum_{w_j} f(w_j)  \ \sum_{w_1, \ldots , w_{j-1}} 
{l_{w_1} \over N} \ {l_{w_2}^{w_1} \over N-1}\ldots
{l_{w_j}^{w_1 \ldots w_{j-1}} \over N-j+1}\cr
&=\sum_{w_j} f(w_j) {l_{w_j} \over N}\cr
&=\bar f
}
$$
where we used lemma 4 that gives 
$$
\sum_{w_{j+1} , \ldots w_n} {l_{w_{j+1}}^{w_1\ldots w_{j}} \over N-j}
\ldots {l_{w_n}^{w_1\ldots w_{n-1}} \over N-n+1} 
=\sum_{w_n} {l_{w_n}^{w_1 \ldots w_j} \over N-j} =1
$$
and 
$$
\sum_{w_1, \ldots , w_{j-1}} 
{l_{w_1} \over N} \ {l_{w_2}^{w_1} \over N-1}\ldots
{l_{w_j}^{w_1 \ldots w_{j-1}} \over N-j+1} ={l_{w_j} \over N}.
$$
}

\stit {c) The true count theorem without induction.}

We present here a direct computation of the expected 
value of the true count when $n$ cards have been withdraw.
It uses the true count theorem for $n=1$ but proves 
the general theorem without induction. 
This is certainly useless, but it prepares the field 
and the notations to prove the standard 
deviation true count formula.

Remember that we denote by $R$ the running count.
The following lemma is immediate by direct computation.

\eno {Lemma 6.}{We have the identity
$$
{R+w_{1}+\ldots +w_{n} \over N-n} -{R\over N} =
{N-1 \over N-n} \left [ \left ( {R+w_{1} \over N-1} -{R\over N} 
\right ) +\ldots + \left ( {R+w_{n} \over N-1} -{R\over N} 
\right ) \right ]
$$
}

Now the expected value of the true count after removing
$n$ cards of the deck is (we use lemma 5 here)
$$
\eqalign {
T_n &= \sum_{w_1, \ldots ,w_n} 
{R+w_{1}+\ldots +w_{n} \over N-n} \ 
{l_{w_1} \over N} \ {l_{w_2}^{w_1} \over N-1}\ldots
{l_{w_n}^{w_1 \ldots w_{n-1}} \over N-n+1} \cr
&= {R \over N}
\ + \sum_{w_1, \ldots ,w_n} \left (
{R+w_{1}+\ldots +w_{n} \over N-n} -{R\over N} \right ) \ 
{l_{w_1} \over N} \ {l_{w_2}^{w_1} \over N-1}\ldots
{l_{w_n}^{w_1 \ldots w_{n-1}} \over N-n+1} \cr
}
$$
Consider the function 
$$
f(w)= {R+w \over N-1} -{R\over N}.
$$
We have (using the true count theorem for $n=1$ (!))
$$
\bar f=0
$$ 
Thus using lemma 6 and lemma 5,  
$$
T_n ={R \over N} + {N-1 \over N-n} n \bar f ={R\over N}.
$$

\stit {d) The true count standard deviation formula.}

We still denote by $f$ the function defined in the 
previous section.
The square of the standard deviation is given by
$$
\s_n^2 =\sum_{w_1, \ldots ,w_n} \left (
{R+w_{1}+\ldots +w_{n} \over N-n} -{R\over N} \right )^2 \ 
{l_{w_1} \over N} \ {l_{w_2}^{w_1} \over N-1}\ldots
{l_{w_n}^{w_1 \ldots w_{n-1}} \over N-n+1} 
$$

Using lemma 6 and developing the square we have 
$$
\s_n^2 = \left ( N-1 \over N-n \right )^2 
\sum_{1\leq j_1 , j_2 \leq n} \ \ \sum_{w_1, \ldots ,w_n}
f(w_{j_1}) f(w_{j_2}) \ 
{l_{w_1} \over N} \ {l_{w_2}^{w_1} \over N-1}\ldots
{l_{w_n}^{w_1 \ldots w_{n-1}} \over N-n+1} 
$$

Immediately from lemma 5 we get

\eno {Lemma 7.}{ For $j_1=j_2=j$, we have 
$$
\sum_{w_1, \ldots ,w_n}
f(w_{j})^2  \ 
{l_{w_1} \over N} \ {l_{w_2}^{w_1} \over N-1}\ldots
{l_{w_n}^{w_1 \ldots w_{n-1}} \over N-n+1}\ =\s_1^2
$$
}

Using next lemma we can finish the computation:
$$
\eqalign {
\s_n^2 &=\left ( N-1 \over N-n \right )^2 \left (n \ \s_1^2 -
n (n+1) {1\over N-1} \ \s_1^2 \right ) \cr
&=\left ( N-1 \over N-n \right )^2 \left ( 1-{n-1 \over N-1} \right )
n \ \s_1^2 \cr
&=\left ( N-1 \over N-n \right ) n \ \s_1^2
}
$$
q.e.d.

\eno {Lemma 8.}{We have for $j_1 \not= j_2$,
$$
\sum_{w_1, \ldots ,w_n}
f(w_{j_1}) f(w_{j_2}) \ 
{l_{w_1} \over N} \ {l_{w_2}^{w_1} \over N-1}\ldots
{l_{w_n}^{w_1 \ldots w_{n-1}} \over N-n+1}\ =-{1\over N-1} \ \s_1^2
$$
}

\dem {Proof.}{We assume $j_1 < j_2$ for example. The proof follows 
the same ideas than the proof of lemma 5. We have using 
lemma 4 several times
$$
\eqalign {
&\sum_{w_1, \ldots ,w_n}
f(w_{j_1}) f(w_{j_2}) \ 
{l_{w_1} \over N} \ {l_{w_2}^{w_1} \over N-1}\ldots
{l_{w_n}^{w_1 \ldots w_{n-1}} \over N-n+1}\cr
&=\sum_{w_1, \ldots ,w_{j_2}}
f(w_{j_1}) f(w_{j_2}) \ 
{l_{w_1} \over N} \ {l_{w_2}^{w_1} \over N-1}\ldots
{l_{w_{j_2}}^{w_1 \ldots w_{j_2-1}} \over N-j_2+1}\cr
&=\sum_{w_{j_2}} f (w_{j_2}) \sum_{w_1 , \ldots , w_{j_2-1}} 
f(w_{j_1}) {l_{w_1} \over N} \ {l_{w_2}^{w_1} \over N-1}\ldots
{l_{w_{j_2}}^{w_1 \ldots w_{j_2-1}} \over N-j_2+1}\cr
&=\sum_{w_{j_2}} f (w_{j_2}) \sum_{w_1 , \ldots , w_{j_1}} 
f(w_{j_1}) \sum_{w_{j_1+1}, \ldots , w_{j_2-1}}
{l_{w_1} \over N} \ {l_{w_2}^{w_1} \over N-1}\ldots
{l_{w_{j_1+1}}^{w_1 \ldots w_{j_1}}\over N-j_1} \ldots 
{l_{w_{j_2}}^{w_1 \ldots w_{j_2-1}} \over N-j_2+1}\cr
&=\sum_{w_{j_2}} f (w_{j_2}) \sum_{w_1 , \ldots , w_{j_1}} 
f(w_{j_1}) {l_{w_1} \over N} \ {l_{w_2}^{w_1} \over N-1}\ldots
{l_{w_{j_1}}^{w_1 \ldots w_{j_1-1}}\over N-j_1+1}\ 
{l_{w_{j_2}}^{w_1 \ldots w_{j_1}} \over N-j_1}\cr
&=\sum_{w_{j_2}} f (w_{j_2}) \sum_{w_{j_1}} f (w_{j_1})
\sum_{w_1 , \ldots , w_{j_1-1}} 
 {l_{w_1} \over N} \ {l_{w_2}^{w_1} \over N-1}\ldots
{l_{w_{j_1}}^{w_1 \ldots w_{j_1-1}} \over N-j_1+1}\ 
{l_{w_{j_2}}^{w_1 \ldots w_{j_1}} \over N-j_1}\cr
&=\sum_{w_{j_1}, w_{j_2}} f(w_{j_1}) f(w_{j_2}) {l_{w_{j_1}}\over N}
\ {l_{w_{j_2}}^{w_{j_1}}\over N-1}\cr
&=-{1\over N-1} \s_1^2
}
$$
In the last computation we use the following lemma 
applied to the function $g(w_1 , w_2) =f(w_1) f(w_2)$
(observe that $\sum_w g(w,w) (l_w /N)=\sum_w f(w)^2 (l_w /N)=\s_1^2$).
}

\eno {Lemma 9.}{Let $g: W\times W \to \dd R$. 
We have 
$$
\sum_{w_1, w_2} g (w_1, w_2) {l_{w_1} \over N} 
{l_{w_2}^{w_1} \over N-1}= {N \over N-1} \ \bar g -{1\over N-1} 
\sum_w g(w,w) {l_w \over N}
$$
where $\bar g$ is the expected value of $g$
$$
\bar g =\sum_{w_1, w_2} g(w_1, w_2) {l_{w_1} \over N} {l_{w_2} \over N}.
$$
}

\dem {Proof.}{We have 
$$
\eqalign {
&\sum_{w_1, w_2} g (w_1, w_2) {l_{w_1} \over N} 
{l_{w_2}^{w_1} \over N-1} \cr
&= \sum_{w_1\not= w_2} 
g(w_1, w_2) {l_{w_1} \over N} 
{l_{w_2}^{w_1} \over N-1} +\sum_w g(w,w) {l_{w} \over N} 
{l_{w} -1 \over N-1}\cr
&= {N\over N-1} \left ( \bar g -\sum_w g(w,w) {l_{w} \over N} 
{l_{w} \over N} \right )
+{N\over N-1} \left (\sum_w g(w,w) {l_{w} \over N} 
{l_{w} -1 \over N-1} \right )\cr
&={N\over N-1} \bar g + {N\over N-1} \left (\sum_w g(w,w) {l_{w} \over N}
\left ( - {1\over N} \right ) \right )\cr
&={N\over N-1} \bar g  -{1\over N-1} \sum_w g(w,w) {l_{w} \over N}\cr
}
$$
}

\null
\vfill
\eject

\tit {3) Long run.}

\stit {a) Kelly for binomial games.}

We review in this section the Kelly criterion in the 
case of an iterated fixed advantage game. We review also
the formulas for the expected value and 
the variance of the exponential rate of 
growth and discuss the implications for Blackjack play.

We assume that we play a repetitive independent 
game where we have 
an advantage $p>1/2$. We expect an exponential growth 
of our initial bankroll $X_0$ if we follow a reasonable
strategy of betting. We assume that there is no minimal 
unit of bet. By homogeneity of the problem, the sharp 
strategy will consists in betting a proportion $f(p)$ of 
the total bankroll. Our bankroll after $n$ rounds of the 
game have been played is
$$
X_n=X_0 \prod_{i=1}^n (1+\eps_i f(p))
$$
where $\eps_i=+1$ if we won the $i$-th hand, and $\eps_i=-1$
if we lost the $i$-th hand.

The exponential rate of growth of the bankroll is
$$
G_n ={1\over n} \log {X_n \over X_0}=
{1\over n} \sum_{i=0}^n \log (1+\eps_i f(p)).
$$
The Kelly criterion maximizes the expected value of the 
exponential rate of growth (the proof is elementary, 
see section b for a proof in the 
more general setting of a "fuzzy advantage").

\eno {Theorem 3.1 (Kelly criterion).}{
The expected value of the exponential rate 
of growth $G_n$ is maximized for 
$$
f(p)=2p-1.
$$
}

Observe that the expected value is
$$
\dd E (G_n)= \dd E (G_1)=p\log (1+f(p))+(1-p) \log (1-f(p)).
$$
Also the random variables 
$$
X_i=\log (1+\eps_i f(p))
$$
are independent, and 
$$
{\rm {\bf Var}} G_n={1\over n}  {\rm {\bf Var}} X .
$$
We compute
$$
\eqalign {
{\rm {\bf Var}}X &=\dd E (X^2)-(\dd E (X))^2 \cr
&=p \left ( \log (1+f(p)) \right )^2 +(1-p) 
\left ( \log (1-f(p)) \right )^2 -\left ( 
p\log (1+f(p))+(1-p) \log (1-f(p)) \right )^2\cr
&=p(1-p) \left (\log \left ({1+f(p) \over 1-f(p)}\right )\right )^2\cr
}
$$
So finally
$$
{\rm {\bf Var}} G_n ={1\over n}p(1-p) \log \left ( {p\over 1-p}\right )^2.
$$

\eno {Proposition 3.2.}{The expected value and the 
standard deviation of $G_n$ are 
$$
\eqalign {
\dd E (G_n) &= \dd E (G_1)=p\log (2p)+(1-p) \log (2-2p) \cr
\s (G_n) &= {1\over \sqrt n} \sqrt {p(1-p)} \log 
\left ( {p\over 1-p}\right )\cr
}
$$
Thus for a slight advantage $p=1/2+\eps$ we have
$$
\eqalign {
\dd E (G_n) &= 2\eps^2 +\cl O (\eps^3)\cr
\s (G_n) &= {2\eps \over \sqrt n} (1+\cl O (\eps^2 ))\cr
}
$$
}

From this data in becomes clear how long is "the long run".
The long run corresponds to play a number of rounds such 
that $\s (G_n)$ is of the order of $\dd E (G_n)$ (or a small fraction
of it). For a smaller number of rounds played
there is a good chance that $G_n$ is negative, i.e. the bankroll
has decreased. In the case of blackjack, $\eps\approx 10^{-2}$. 
Thus the "long run" corresponds to have played of the order 
of $10 000$ favorable hands, thus at least $20 000$ hands. At
$50$ hands per hour this adds to $4 000$ hours of play. To play less 
hands means that we are gambling. This explains that even at 
the non-professional level the team play makes perfect sense in 
order to get into the long run quicker (divide the time 
by the number of members of the team).

The situation one faces when playing blackjack is different 
than a binomial game with a fixed advantage. First the advantage
is not the same at different hands. This involves minor modifications
of the above computations. Second, due to the standard deviation
effect, the advantage the player has at the moment of playing is 
a random variable at the moment of betting. So he is betting 
according to a fuzzy advantage. We study this situation in the
next sections. Kelly criterion is still sharp. The main effect 
of the fuzzy advantage is to introduce a supplementary 
term in the standard deviations of $G_n$.

\stit {b) Kelly for fuzzy advantage.}

We consider the situation one faces playing blackjack. At 
each round one decides the amount to bet in function of the 
true count.

We assume that we play multiple independent
rounds of the same game. At each round the advantage we 
have is a random variable $p$, $0 < p <1$, with a known distribution
$D_x(p)$ and expected value $p_0(x) >1/2$ from a family of 
distributions $(D_x)$. The choice of the distribution $D_x$ at 
each round is random with distribution $\rho (x)$, $x\in [0,1]$ ($x$
is a mere index, we may just adjust it to have a uniform 
distribution $\rho (x)=1$). 

We want to maximize the expected value of 
the exponential rate of growth
of our bankroll $X_0$. The optimal 
amount to bet
is a proportion $0\leq f(D_x) \leq 1$ 
of the bankroll depending only on the distribution
$D_x$. For a sharp strategy we have $f(D_x)=0$ when 
$p_0 (x) \leq 1/2$ (that is, the game is not favorable). 
The quantity to maximize is 
$$
\eqalign {
&\int_0^1 \int_0^1 \left ( p\log (1+f(D_x)) +(1-p) \log (1-f(D_x))\right )
\ D_x(p) \ dp \ \rho (x) \ dx \cr
&= \int_0^1 \left ( p_0 (x) \log (1+f(D_x)) +(1-p_0 (x)) 
\log (1-f(D_x))\right )
\ \rho (x) \ dx \cr
}
$$

\eno {Theorem 3.3 (Kelly criterion).}{ 
The optimal strategy is obtained for
$$
f_K (D_x) =2 p_0 (x) -1
$$
for $1/2\leq p_0(x) \leq 1$, and $f_K (D_x) =0$ for 
$0\leq p_0(x)\leq 1/2$.
}

\dem {Proof.}{We want to maximize the functional 
$$
\eqalign {
\cl G (f) &= \int_0^1 \left ( p_0 (x) \log (1+f(D_x)) +(1-p_0 (x)) 
\log (1-f(D_x))\right )
\ \rho (x) \ dx \cr
&=\int_{\{x; p_0 (x) >1/2 \}} \left ( p_0 (x) \log (1+f(D_x)) +(1-p_0 (x)) 
\log (1-f(D_x))\right )
\ \rho (x) \ dx \cr}
$$
where we assume that $f(D_x)=0$ for $0<p_0 (x) \leq 1/2$.
If $f_0$ is an extremum then for any perturbation $h$ such that
$f_0+\eps h$ is an allowable strategy ($f_0 +\eps h >0$), then if 
we consider 
$$
g(\eps) =\cl G (f_0 +\eps h)
$$
we must have 
$$
g'(0)=0.
$$
A direct computation gives
$$
g'(\eps )=\int_{\{x; p_0 (x)>1/2 \} } 
{h(D_x) \over 1-\left ( f_0(D_x)+\eps h(D_x) \right )^2}
[(2p_0 (x)-1)-f_0(D_x) -\eps h(D_x) ] \ \rho (x) \ dx>
$$
So
$$
g'(0) =\int_{\{x; p_0 (x)>1/2 \} } 
{h(D_x) \over 1-f_0 (D_x)^2} \ (2 p_0 (x)-1 -f_0 (x))
\ \rho (x) \ dx \ .
$$
Thus when $p_0(x) >1/2$ we must have $f_0 (x)=2p_0 (x)-1 =f_K (x)$.
Also by direct computation we have
$$
g''(\eps )= \int_{\{x; p_0 (x)>1/2 \} } {- \ h^2 \over (1-(f_0+\eps h)^2)^2}
\left [ (f_0+\eps h)^2 +1 -2 (2p_0(x)-1)(f_0+\eps h) \right ]\ \rho (x) \ dx
$$
Now since $f_0+\eps h >0$ and $2 p_0 (x)-1 >0$ in the range of 
integration, we have
$$
g''(\eps ) \geq \int_{\{x; p_0 (x)>1/2 \} } {h^2 \over (1-(f_0+\eps h)^2)^2}
(f_0+\eps h -1)^2 \rho (x) \ dx >0
$$
Therefore $\eps \mapsto g(\eps )$ is concave, $g''(0) < 0$, and 
$f_0=f_K$ is a maximum of the functional.}

\stit {c) Long run.}

We assume here that the player uses a Kelly betting strategy
for a repetitive game with fuzzy advantage.
For a given advantage distribution $D_x$ with $p_0 (x) >1/2$, 
the expected exponential
rate of growth after playing $n$ favorable hands with distributions
$x_1, x_2, \ldots$ is 
$$
{\bf E}(G_n)={1\over n} \sum_{i=1}^n \int (p\log (1+f(p_0 (x_i)))+
(1-p) \log (1-f(p_0(x_i))))\ D_{x_i} (p) \ dp.
$$
Thus the expected exponential rate of growth is
$$
\eqalign {
{\bf E} (G_n) &=\int \int (p\log (1+f(p_0 (x)))+
(1-p) \log (1-f(p_0(x))))\ D_{x} (p) \ dp \rho (x) \ dx \cr
&=\int (p_0(x) \log (1+f(p_0 (x)))+(1-p) \log (1-f(p_0(x))))
\ \rho (x) \ dx \cr
&= \int {\bf E}_x (G_1) \ \rho (x) \  dx \cr
&= {\bf E} (G_1)\cr
}
$$

Observe also that by independence of the hands we have
$$
{\bf Var } (G_n) ={1\over n} {\bf Var} (G_1).
$$

We extimate the standard deviation of $G_n$. For this we make
the assumption that the distribution $D_x(p)$ is a distribution
$D(p)$ with a small "noise", that is
$$
D_x (p) =D(p) + Z_x \ d(p)
$$
where $Z_x$ is a random variable with $0$ expectation. 
We have
$$
\eqalign {
\int Z_x \rho (x) dx =0 \cr
\int d(p) dp =0 \cr
}
$$
the second equation coming from $\int D_x (p)\  dp=1$.
According to the previous notation, the expected advantage
is
$$
p_0 (x) = \int p D_x (p) \ dp
$$
thus 
$$
p_0 (x) =p_0 +A Z_x 
$$
where $p_0 =\int p D(p) \ dp$ and we denote
$$
A=\int p d(p) \ dp.
$$
Observe also that the quantity $A$ is directly related to the 
variance of $p_0 (x)$ by
$$
{\bf Var }( p_0 (x)) = \int p_0 (x)^2 \rho (x) dx -p_0^2 =
A^2 \int Z_x^2 \rho (x) dx \ .
$$

Our purpose now is to compute the first order perturbation of 
the standard deviation of $G_n$ introduced by the noise $Z_x$.

We have to the second order on $Z_x$ (we do not care about the 
coefficient of the first order),
$$
\eqalign {
{\bf E}_x (G_1) &= (p_0 +A Z_x) \log (2p_0 +2A Z_x) +
(1-p_0-A Z_x) \log (2-2p_0 -2A Z_x) \cr
&= (p_0 +A Z_x)  (\log (2p_0) +\log (1+A Z_x /p_0)) +\cr
& +(1-p_0-A Z_x) (\log (2-2p_0)+\log (1 -AZ_x /(1-p_0))) \cr
&= p_0 \log (2p_0) +(1-p_0) \log (2-2p_0) + * AZ_x +
{1\over 2 p_0 (1-p_0)} A^2 Z_x^2 +\ldots \cr
}
$$
Thus taking expected values with respect to $x$,
$$
{\bf E} (G_1) =p_0 \log (2p_0) +(1-p_0) \log (2-2p_0)
+{1\over 2} {{\bf Var} (p_0 (x)) \over p_0 (1-p_0)} +\ldots 
$$
And, on the first order on ${\bf Var}  (p_0 (x))$, we have 
$$
{\bf E} (G_1)^2 = (p_0 \log (2p_0) +(1-p_0) \log (2-2p_0))^2
+(p_0 \log (2p_0) +(1-p_0) \log (2-2p_0))
{{\bf Var} \ (p_0 (x)) \over p_0 (1-p_0)} +\ldots 
$$

Now we compute the expansion of
$$
{\bf E} (G_1^2) =\int {\bf E}_x (G_1^2) \rho (x) \ dx \ .
$$
We have (after a long computation)
$$
\eqalign {
{\bf E}_x (G_1^2) &=p_0 (x) (\log (2p_0 (x)))^2 +
(1-p_0 (x)) (\log (2-2p_0(x)))^2 \cr
&=\ldots \cr
&= p_0 (\log (2p_0))^2 +
(1-p_0) (\log (2-2p_0))^2 + *\  AZ_x + \cr
&+ {1\over p_0 (1-p_0)} \left ( 1+ (1-p_0) \log (2p_0) +p_0 
\log (2-2p_0) \right ) A^2 Z_x^2 +\ldots  \cr
}
$$
Thus 
$$
\eqalign {
{\bf E} (G_1^2) &=
 p_0 (\log (2p_0))^2 +
(1-p_0) (\log (2-2p_0(x)))^2 + \cr
& +{1\over p_0 (1-p_0)} \left ( 1+ (1-p_0) \log (2p_0) +p_0 
\log (2-2p_0) \right ) {\bf Var} (p_0(x)) +\ldots \cr
}
$$

Finally, putting together the previous formulas, we have

\eno {Theorem 3.4.}{The variance of the exponential
rate of growth in a Kelly betting strategy with fuzzy advantage
is, at the first order in the noise,
$$
{\bf Var} (G_1) =p_0 (1-p_0) \left ( \log \left ( {p_0 \over 1-p_0}
\right ) \right )^2 +{1\over p_0 (1-p_0)}\left ( 1-(2p_0-1)
\log \left ({p_0 \over 1-p_0}\right ) \right ) {\bf Var} (p_0 (x))
+\ldots 
$$
}

As an appication to blackjack, where we have $p_0 =1/2+\eps$
with $\eps \approx 10^{-2}$, and ${\bf Var} (p_0 (x)) =\eps \s_{BET}$
(where $\s_{BET}$ is the standard deviation of the true count),
$$
{\bf Var} (G_1) =4 (1+\s_{BET}^2 )\eps^2 +\ldots 
$$
thus 
$$
\s (G_n) ={2 \sqrt {1+\s_{BET}^2} \over \sqrt n} \eps +\ldots
$$
Thus the long run is increased by a factor $1+\s_{BET}^2$
caused by the standard deviation effect. 
In order to carry out precise estimates on the "long run" 
we state a precise definition. Note that the factor
$2$ is somewhat arbitrary. We may want to choose another 
value in order to have a bound on the probability of loosing
after achieving the long run $N$ using Tchebichev's inequality.

\eno {Definition 3.5 (Long run).}{The long run $N$ is the 
minimal number of hands in order to have
$$
{{\bf E } (G_N) \over \s (G_N)} \geq 2.
$$
}

Thus from the above computation we have (disregarding integer 
parts),
$$
N=2^2 (1+\s_{BET}^2) \eps^{-2}\approx (1+\s_{BET}^2) .  40000
$$
Typically the difference between $\s_{BET}^2$ between first and 
third base could be of the order of $2 \%$ which makes the 
long run about $2 \% $ longer, or about $800$ more favorables 
hands to be played, thus about $2000$ more hands to be played
(that is at least $40$ more hours of play).

\null
\vfill
\eject

\tit {4) Practical gambling.}

\stit {a) Comparison of first and third base and head-on play.}

\stit {a.1) Full table.}

The number of cards removed from the deck 
between the betting decision and 
the play decision of the first base is exactly $8\times 2 =16$.
The average number of cards for the third base can be computed 
considering that the average number of cards in a 
hand of players {\bf playing
basic strategy} is 2.6 (this number is well known). 
Thus in average the number of cards 
removed from the deck between the betting decision and the 
play decision is 19.6. Denoting by $\Sigma (1)$ and 
$\Sigma (2)$ the corresponding standard deviations of the 
true count, and considering the approximation $N >>n$, we have
$$
{\Sigma (2) \over \Sigma (1)} \approx {4.42 \over 4}=1.10
$$
Thus $\Sigma (2)$ is in average about $10\% $ higher than 
$\Sigma (1)$.

\stit {In company of ploppys.}

A recurrent topic in the blackjack literature is about 
is the influence of inaccurate players (ploppys) 
in the same table.
It is an easy escape gate to blame others for your losses.
The systematic answer that one founds all 
over the blackjack literature is that this is pure 
superstition. We believe also that there is mostly exaggeration
in these complaints, but as we explain next, there is some mathematical 
foundation based on the standard deviation effect.

\medskip

In the first part of this section, when we were comparing the 
standard deviation of the first and third base, in order
to estimate the number of cards played between first and 
third base, we made a major assumption: Other players at 
the table are playing basic strategy, thus the average number
of cards in a hand is 2.6. Regular players report that very 
few players know basic strategy (less than $5 \% $ we were told)
and among these maybe about $1/4$ of them are more or less 
skilled card counters. Playing basic strategy, the casino has 
an edge of about $1 \% $ ($0.5 \% $ with best rules), but the 
actual win rate of blackjack in a casino is estimated to 
be between $3 \% $ and $7 \%$ (depending of the sources, 
see [Gr2] p.137). This shows 
how inaccurately the average player plays blackjack. The most 
common mistake is to think that the game consists in approaching
as accurately as possible a total sum of $21$. This is totally 
false. Doing so one risks being busted, thus loosing the whole
bet without getting the chance of letting the dealer bust. 
The real goal is to beat the dealer. And the best way 
to achieve this is to let him bust. Also ploppys like to split 
tens for example. This is something at which strategy players 
look with horror because the count has to be sky high to 
justify this play (nevertheless 
it can happen, see [Us] p.17 where J. Uston explain how he,  
on one occasion, he did split
eight times tens in one hand...with an initial bet of \$ 1000). 
This mistakes make that 
the average player uses more cards in his hands than the 
ideal basic strategy player. We do not 
have precise statistics on this, but one can guess that having 
a gang of ploppys spliting tens to your right can make  
that if you sit in third base the average number of cards played between
the betting and the play decision could go as high as $25$.
In such a table the standard deviation of third base will 
be $25 \% $ higher than the one for first base. Thus this will 
induce longer long run and larger fluctuations, 
and you will be right to blame the 
ploppys.

\stit {a.2) Head-on play.}

The differences between the standard deviation effect when 
one plays head-on (that is alone with the dealer)
and one plays with a full table is even more important.
The exact number of cards played in between decisions 
for a head-on play is $4$. Thus if you sit in first 
base at a full table your standard deviation will be twice
as important (100 \% more). And if you sit in third base
this figure goes up to $121 \% $ more.

\stit {b) Absolute magnitude of standard deviation.}

The standard deviation true count formula allows us to compute
$\s_n$ from $\s_1$. We compute now the standard deviation 
$\s_1$. It not only depends on $N$, and $R$, and the system 
count but also in the actual composition of the remaining deck.

\eno {Proposition 1.}{We have 
$$
\s_1 ={1\over N-1} \sqrt { \sum_w w^2 \ {l_w \over N} -\left (
{R\over N} \right )^2 }.
$$
}

\dem {Proof.}{We have 
$$
\eqalign {
\s_1^2 &= \sum_w \left ( {R+w \over N-1} -{R\over N} \right )^2\ 
{l_w \over N} \cr
&={1\over N^2 (N-1)^2} \sum_w (N w+R)^2 {l_w \over N} \cr
&= {1\over N^2 (N-1)^2} \left [ N \left ( \sum_w w^2 l_w \right )
+2R \sum_w w l_w +R^2 \right ] \cr
&={1\over N^2 (N-1)^2} \left [ N \left ( \sum_w w^2 l_w \right )
-R^2 \right ]\cr
&={1\over (N-1)^2} \left [ N \left ( \sum_w w^2 l_w \right )
-\left ({R \over N} \right )^2 \right ]\cr
}
$$
}

Let $\cl S$ denote the balanced counting system used. 
A main characteristic 
of the counting system 
is the standard deviation $\Sigma_0 (\cl S)$ of the weights.
Since the system is balanced we have
$$
\Sigma_0 = \Sigma_0 (\cl S) =\sum_w w^2 {S_w \over 52}
$$
where $S_w$ is the total number of cards of weight $w$ in a complete 
deck. There follows the standard deviation 
of different systems round up to the third digit (for a
description the different systems with their list of weights 
the reader can consult [Sc] p.62).

\bigskip

\item {$\bullet$} Uston ace-five count: \hfill $\Sigma_0 =0.392 $ 
\hskip 2cm \null 

\item {$\bullet$} Hi-Opt I: \hfill $\Sigma_0 = 0.784$ \hskip 2cm \null 

\item {$\bullet$} C-R Point count: \hfill $\Sigma_0 =0.855 $ \hskip 2cm \null 

\item {$\bullet$} Canfield expert count: \hfill $\Sigma_0 =0.877 $
\hskip 2cm \null 

\item {$\bullet$} Hi-Lo: \hfill $\Sigma_0 =0.877$ \hskip 2cm \null 

\item {$\bullet$} Uston advanced plus-minus count: \hfill 
$\Sigma_0 =0.877 $ \hskip 2cm \null

\item {$\bullet$} Halves count: \hfill $\Sigma_0 =0.920 $ \hskip 2cm \null 

\item {$\bullet$} The systematic count: \hfill $\Sigma_0 =0.961 $ 
\hskip 2cm \null 

\item {$\bullet$} Hi-Opt II: \hfill $\Sigma_0 =1.468 $ \hskip 2cm \null 

\item {$\bullet$} Canfield master count: \hfill $\Sigma_0 =1.569 $
\hskip 2cm \null 

\item {$\bullet$} Zen count: \hfill $\Sigma_0 =1.569$ \hskip 2cm \null

\item {$\bullet$} Uston advanced point count: \hfill 
$\Sigma_0 =1.687$ \hskip 2cm \null

\item {$\bullet$} The Revere point count: \hfill $\Sigma_0 = 1.710$
\hskip 2cm \null

\item {$\bullet$} Uston SS count: \hfill $\Sigma_0 =1.797 $ \hskip 2cm \null

\item {$\bullet$} The Revere advanced point count: \hfill 
$\Sigma_0 =2.449$ \hskip 2cm \null

\item {$\bullet$} Griffin seven count: \hfill 
$\Sigma_0 =3.234 $ \hskip 2cm \null

\item {$\bullet$} Thorp ultimate count: \hfill  $\Sigma_0 =5.798 $
\hskip 2cm \null

\medskip

We observe, as expected, that $\Sigma_0$ is larger for 
higher level (i.e. more precise) count systems.

From the formula in proposition 1 we can get a good 
approximation of $\s_1$ in real gambling situations.
Observe that for shoe games $N$ is always large due to 
non complete penetration (typically for shoe games we 
always have $N\geq 50$). Thus the square of the standard 
deviation in the remaining deck, that is 
$$
\sum_w w^2 {l_w \over N}  \ ,
$$
is well approximated by $\Sigma_0^2$. Thus the square
of the true count $(R/N)^2$ can be neglected (remember 
that the units we use to compute the number of cards is 
card units and not deck units as is done in practice).
We conclude that the quantity
$$
\s_1' ={\Sigma_0 \over N}
$$
is a good approximation of $\s_1$ provided that the 
shoe is not totally depleted (note that in the approximation
we replace $N-1$ by $N$ which is acceptable).

In practice one measures $N$ and the true count 
in deck units, then using 
the previous formula one will get $\s'_1$ in deck units
that is the standard deviation of the true count in deck 
units which is the quantity that is used to determine the 
advantage in percentage.
We can also approximate $\s_n$:
$$
\s_n =\sqrt {N-1 \over N-n} \ \sqrt n \ \s_1 
\approx \sqrt n  \ {\Sigma_0 \over N}
$$
where we did approximate $\sqrt {(N-1) /(N-n)}$ by $1$, which 
is quite accurate (except maybe if just one deck remains, then 
the value of $\s_n$ will be larger than the approximation).

Using this approximation one can compute the following 
tables for standard deviation related to the bet decision 
($\s_{BET}$) and the one for the play decision ($\s_{PLAY}$),
for different players and 
different types of games.

\bigskip

$$
\vbox{
\offinterlineskip
\halign{
\strut \vrule $#$ &\vrule $#$ &\vrule $#$ &\vrule $#$ \vrule \cr
\noalign{\hrule}
{\hbox { Position }} & \ 1 & \ 4 & \ 7\cr
\noalign{\hrule}
\ \sigma_{BET} \  & \ 0.877 & \ 0.925 & \ 0.971 \cr
\noalign{\hrule}
\ \sigma_{PLAY} \ & \ 0.449 & \ 0.340 & \ 0.085 \cr
\noalign{\hrule}
}
}
$$

\medskip

{\centerline {\bf Hi-Lo, Canfield expert, 
Uston Advanced plus-minus ($\Sigma_0 = 0.877$)}}

{\centerline {\bf $8$-Deck shoe, $50 \% $ shoe played}}

\bigskip

$$
\vbox{
\offinterlineskip
\halign{
\strut \vrule $#$ &\vrule $#$ &\vrule $#$ &\vrule $#$ \vrule \cr
\noalign{\hrule}
{\hbox { Position }} & \ 1 & \ 4 & \ 7\cr
\noalign{\hrule}
\ \sigma_{BET} \  & \ 1.316 & \ 1.388 & \ 1.457 \cr
\noalign{\hrule}
\ \sigma_{PLAY} \ & \ 0.674 & \ 0.510 & \ 0.128 \cr
\noalign{\hrule}
}
}
$$

\medskip

{\centerline {\bf Hi-Lo, Canfield expert, 
Uston Advanced plus-minus ($\Sigma_0 = 0.877$)}}

{\centerline {\bf $8$-Deck shoe, $66.6 \% $ shoe played}}

\bigskip

$$
\vbox{
\offinterlineskip
\halign{
\strut \vrule $#$ &\vrule $#$ &\vrule $#$ &\vrule $#$ \vrule \cr
\noalign{\hrule}
{\hbox { Position }} & \ 1 & \ 4 & \ 7\cr
\noalign{\hrule}
\ \sigma_{BET} \  & \ 1.754 & \ 1.850 & \ 1.942 \cr
\noalign{\hrule}
\ \sigma_{PLAY} \ & \ 0.898 & \ 0.680 & \ 0.170 \cr
\noalign{\hrule}
}
}
$$

\medskip

{\centerline {\bf Hi-Lo, Canfield expert, 
Uston Advanced plus-minus ($\Sigma_0 = 0.877$)}}

{\centerline {\bf $8$-Deck shoe, $75 \% $ shoe played}}

\bigskip

$$
\vbox{
\offinterlineskip
\halign{
\strut \vrule $#$ &\vrule $#$ &\vrule $#$ &\vrule $#$ \vrule \cr
\noalign{\hrule}
{\hbox { Position }} & \ 1 & \ 4 & \ 7\cr
\noalign{\hrule}
\ \sigma_{BET} \  & \ 5.780 & \ 6.096 & \ 6.400 \cr
\noalign{\hrule}
\ \sigma_{PLAY} \ & \ 2.959 & \ 2.241 & \ 0.560 \cr
\noalign{\hrule}
}
}
$$

\medskip

{\centerline {\bf Thorp ultimate ($\Sigma_0 = 5.798$)}}

{\centerline {\bf $8$-Deck shoe, $50 \% $ shoe played}}

\bigskip

$$
\vbox{
\offinterlineskip
\halign{
\strut \vrule $#$ &\vrule $#$ &\vrule $#$ &\vrule $#$ \vrule \cr
\noalign{\hrule}
{\hbox { Position }} & \ 1 & \ 4 & \ 7\cr
\noalign{\hrule}
\ \sigma_{BET} \  & \ 8.670 & \ 9.144 & \ 9.600 \cr
\noalign{\hrule}
\ \sigma_{PLAY} \ & \ 4.439 & \ 3.362 & \ 0.840 \cr
\noalign{\hrule}
}
}
$$

\medskip

{\centerline {\bf Thorp ultimate ($\Sigma_0 = 5.798$)}}

{\centerline {\bf $8$-Deck shoe, $66.6 \% $ shoe played}}

\bigskip

$$
\vbox{
\offinterlineskip
\halign{
\strut \vrule $#$ &\vrule $#$ &\vrule $#$ &\vrule $#$ \vrule \cr
\noalign{\hrule}
{\hbox { Position }} & \ 1 & \ 4 & \ 7\cr
\noalign{\hrule}
\ \sigma_{BET} \  & \ 11.56 & \ 12.192 & \ 12.8 \cr
\noalign{\hrule}
\ \sigma_{PLAY} \ & \ 5.918 & \ 4.482 & \ 1.120 \cr
\noalign{\hrule}
}
}
$$

\medskip

{\centerline {\bf Thorp ultimate ($\Sigma_0 = 5.798$)}}

{\centerline {\bf $8$-Deck shoe, $75 \% $ shoe played}}

\bigskip

A common feature is that the fluctuations become much larger
at the end of the deck. And there is where most of the action 
takes place because it is at this moment than the true count 
uses to take its largest values ! This picture is scaring 
for the card counter and explains why blackjack is a game
with high fluctuations for card counters. Let's see what 
Stanford Wong and others have to say on this subject in his reference 
book [Wo] (p.199):

\medskip

{\it {\bf Comments on risks.} Peter Giles said: "The late Ken Uston 
was once quoted in the "Review Journal" as saying, "It's really
tough to make a living at blackjack. The fluctuations will really 
wipe out the average guy. If I had to play  by myself (instead of 
on a team), I probably wouldn't be in it now". You can quote 
me (S. Wong) as saying the same (...) I am still trying to 
determine how many units one is safe with. What is recommended 
in books is, in my opinion, too risky. This is one area in 
which it is hard to trust mathematics.}

\medskip

Probably not all the mathematics can be found in the classical
litterature. It is surprising to note that there is no complete
treatment of the long run in the blackjack litterature as the one
we carry out in section 3.

To limit risks some card counters
use to divide their true count by the number of players at the 
table (as Sonia, a professional gambler, does [So]), 
thus limiting the action, the fluctuation, 
... but also the advantage. 
The standard deviation tends 
to infinite when the deck depletes completely. We can see this 
formally, making $N \to 0$ in the formulas (in a 
continuous model $N-1$ is replaced by $N$). 
Of course $N$ takes integer values
but we can compute for example that if only $1/4$ of a deck remains
and we assume that $\Sigma_0 \approx 1$ for the count system used,
then 
$$
\s_1 \approx 4
$$
thus in the best case (playing head-on),
$$
\Sigma (1)\approx 8.
$$
And in the case of a full table, third base placement,
$$
\Sigma (2) \approx 17.7,
$$
a true count of the order of 17 is certainly like $+\infty$.

When one sees these figures one wonders how it is possible that 
removing just one card from the deck could have such an effect
in the standard deviation. Let's compute exactly $\s_1$ in one 
particular situation. We assume that in a head-on game
$13$ cards remain (that is $1/4$ of a deck). 
We assume that we use Hi-Lo and the the composing of the 
remaining cards is $5$ high cards, $5$ low cards and $3$ medium
cards (thus $R=0$ at this point). The standard deviation of the 
true count after removing just one card is (in deck units)
$$
\s_1=52. \sqrt { \left ( {1\over 12}\right )^2 \ {5 \over 13}+ 
\left ( {1\over 12}\right )^2 \ {5 \over 13} +0 } =52. 
\sqrt {5 \over 936} \approx 3.80
$$
And if still one is not convinced, just look at what happens 
when we reveal one card, and this card is for example a low 
card. The new true count (in deck units) will be
$$
52.{R\over N} ={52 \over 12}\approx 4.33
$$

That means that in such situations the true count indicator
is meaningless and the standard deviation effect takes over. 
In an evil scenario casinos could exploit 
that weakness of card counters. They could offer and advertise 
very good penetration in order to get tables crowded with card counters.
Then deal the shoes almost to bottom to exploit the huge standard 
deviation effect. Cardcounters will then experience too large 
fluctuations that will bring them into an important risk of
ruin. Thus the casino will wipe out bankrolls 
of some unlucky card counters. It is unlikely that the casino will
make any money in the operation, due to the increased 
benefits of the survival
card counters.

\null
\vfill
\eject

{\centerline {\bf {Bibliography}}}

\bigskip
\bigskip


\ref{[Br]}{L. BREIMAN}{Optimal gambling systems for favorable
games}
{Fourth Berkeley Symposium on Probability and Statistics, {\bf I},
p. 65-78, 1961.}

\ref{[Gr1]}{P.A. GRIFFIN}{The theory of blackjack}
{Huntington Press, Las Vegas, 1999.}

\ref{[Gr2]}{P.A. GRIFFIN}{Extra stuff}
{Huntington Press, Las Vegas, 1991.}

\ref{[HC]}{L. HUMBLE, C. COOPER}{The world's greatest blackjack book}
{Doubleday, New York,  1980.}

\ref{[JM]}{A. JALIB M'HALL}{The true count theorem}
{Message posted in rec.gambling.blackjack, 7/30/1996.}


\ref{[RT]}{L.M. ROTANDO, E.O. THORP}{The Kelly criterion and the 
stock market}
{Amer. Math. Monthly, {\bf 99}, p.922-931, 1992.}

\ref{[Sc]}{F. SCOBLETE}{Best blackjack}
{Bonus Books, Chicago, 1996.}

\ref{[So]}{SONIA}{Personal communication}
{1/2000.}

\ref{[Th]}{E. O. THORP}{Beat the dealer}
{Vintage books, New York, 1962.}

\ref{[Wo]}{S. WONG}{Professional blackjack}
{Pi Yee Press, 1994.}

\end